\documentclass{article}
\usepackage[utf8]{inputenc}
\usepackage{authblk}
\usepackage[colorlinks]{hyperref}

\title{On a General Method for Resolving Integrals of Multiple Spherical Bessel Functions Against Power Laws into Distributions}

\author[1*]{Kiersten Meigs}
\author[1,2]{Zachary Slepian}
\affil[1]{\footnotesize Department of Astronomy, University of Florida, 211 Bryant Space Science Center, Gainesville, FL 32611, USA}
\affil[2]{\footnotesize Physics Division, Lawrence Berkeley National Laboratory, 1 Cyclotron Road, Berkeley, CA 94709, USA}
\affil[*]{Electronic Address: \href{mailto: meigsk@ufl.edu}{meigsk@ufl.edu}}

\usepackage[a4paper, total={6in, 8in}]{geometry}
\usepackage{amsmath}
\usepackage{amssymb}

\usepackage[dvipsnames]{xcolor}

\usepackage[numbers]{natbib}
\usepackage{graphicx}
\usepackage{subcaption}

\newcommand{\dD}{\delta_{\rm D}^{[1]}}

\usepackage{color}
\usepackage{enumitem}
\def\beq{\begin{eqnarray}}
\def\eeq{\end{eqnarray}}

\definecolor{darkgreen}{RGB}{0,120,0}

\begin{document}

\maketitle

\begin{abstract}
We here present a method of performing integrals of
products of spherical Bessel functions (SBFs) weighted by a power-law. Our method, which begins with double-SBF integrals, exploits a differential operator $\hat{D}$ defined via Bessel's differential equation. Application of this operator raises the power-law in steps of two. We also here display a suitable base integral expression to which this operator can be applied for both even and odd cases. We test our method by showing that it reproduces previously-known solutions. Importantly, it also goes beyond them, offering solutions in terms of singular distributions, Heaviside functions, and Gauss's hypergeometric,$\;_2{\rm F}_1$ for \textit{all} double-SBF integrals with positive semi-definite integer power-law weight. We then show how our method for double-SBF integrals enables evaluating \textit{arbitrary} triple-SBF overlap integrals, going beyond the cases currently in the literature. This in turn enables reduction of arbitrary quadruple, quintuple, and sextuple-SBF integrals and beyond into tractable forms.
\end{abstract}

\section{Introduction}
\label{sec:intro}
\subsection{Motivation}
Integrals of spherical Bessel functions (SBFs) are ubiquitous in cosmology. Integrals against $j_0$ are used to transform the power spectrum into the correlation function, or the bispectrum into the 3-Point Correlation Function (\textit{e.g.} \cite{model_3}, \cite{umeh}, \cite{aviles21}), and integrals against $j_1$ and $j_2$ respectively enter models involving either a vector or a tensor dependence of the model (\textit{e.g.} \cite{rv}). They also appear in exact consideration of fully spherical redshift-space distortions (\textit{e.g.} \cite{hamilton}, \cite{castorina}). Integrals of two SBFs enter both the covariance matrix of the correlation function (\textit{e.g.} \cite{xu}) and the computation of the CMB lensing angular power spectrum (often done using the Limber approximation). Double SBF integrals, and even triple SBF integrals, also enter computation of the galaxy 3PCF and 4PCF covariance matrices (\cite{3pt_alg}, \cite{aniso_3}, \cite{jiamin_covar}) and indeed it can be shown that the Gaussian Random Field approximation of the covariance for any N-Point Correlation Function is expressible in terms of double and triple SBF integrals against the power spectrum \cite{jiamin_covar}. Famously, separating the double-SBF integrals against the source function into a ``geometric'' and ``numerical'' part is at the heart of the acceleration of the CMBFAST algorithm for CMB anisotropies \cite{uros}. Double and triple SBF integrals also enter the computation of loop corrections to the power spectrum \cite{schmitt_1, schmitt_2, pt_decoup} in perturbation theory and the redshift-space power spectrum in Effective Field Theory of Large-Scale Structure \cite{fonseca}. Recently, there has been a movement to expand the power spectrum into complex power laws and then perform the SBF-against-power-law integrals analytically and sum them with the appropriate expansion-coefficient weights just once at the end. This approach has been used both to accelerate computing perturbation theory \cite{simonovic} and angular statistics models (both primary CMB \cite{schoneberg}, lensing \cite{beyond_limber} and projected \cite{assassi}, \cite{2fast}) and also in computing the relevant covariance matrices \cite{fang}. SBF overlap integrals even appear in such diverse cosmological contexts as a toy model for the effect of a survey boundary \cite{hand}, analytic treatment of the development of the Baryon Acoustic Oscillation feature in the matter transfer function \cite{se_bao_analyt}, and faster computation of the root-mean-square mass fluctuations on a given smoothing scale \cite{krol}.

With these many reasons for interest in such integrals, it is worthwhile to expand our arsenal of techniques for performing them. In particular, a lacuna in the literature has been resolution of singular versions of such integrals, \textit{i.e.} ones where the power-law weight is sufficient to make the integral diverge in the UV. Of course, a separate matter of interest might be IR-divergent such integrals, but we leave this unconsidered here. This lacuna has developed because, while Dirac delta functions are of great use to physicists, mathematicians tend to avoid them as defective functions. Yet, the types of singular integrals we treat in the present work can arise in cosmological Perturbation Theory (PT) (\textit{e.g.} Hou \& Slepian in prep.). These integrals are also of intrinsic interest.

\subsection{Prior work: convergent double and triple sBF integrals}
Prior work on such integrals has been extensive and we summarize it here. We organize our discussion by similar integrals; we first focus on the most recent results and then proceed chronologically. Though little has been done on indefinite integrals and they are not our concern here, we note for completeness the work in this regard of \cite{indef_bloom} and \cite{indef_teboho}. 

Regarding definite integrals, the series of papers by Mehrem has been of significant use and generalizes beyond early results \textit{e.g.} in \cite{bailey} or in Watson's Treatise \cite{wat}. Mehrem shows that a triple-SBF integral satisfying certain conditions on the orders and the power can be written as a sum over angular momentum symbols weighting Legendre polynomials in, essentially, the cosine of the angle enclosed by the triangle formed by the three free arguments of the integrand (\textit{e.g.} $r_1, r_2$ and $r_3$ in $\int k^2 dk\; j_{\ell_1}(k r_1) j_{\ell_2}(k r_2) j_{\ell_3}(k r_3)$) \cite{mehrem_pwe}, \cite{mehrem_3_gen}. Mehrem also included in double-SBF integrals \cite{mehrem_2_91} additional factors in the integrand beyond power-laws, such as an exponential and polynomial \cite{mehrem_2_exp_poly} and a Gaussian \cite{mehrem_2_gauss}. Importantly, Mehrem also displayed useful relationships between the plane-wave expansion and such integrals \cite{mehrem_pwe}. Though his work exploits that a Dirac delta function can be expressed as a Fourier Transform of unity, and hence plane-wave expansions can be invoked, his end results focus on only convergent integrals.

Much earlier work exploiting the plane-wave expansion is \cite{jackson1972integrals}, which obtains integrals of three spherical Bessel functions with $k^2$ weight and subject to the condition that the sum of the orders is even and that the orders satisfy triangular inequalities; if these two conditions hold, then the integrals vanish unless the three free arguments also can form a triangle. \cite{gervois85} obtains solutions for integrals of two Bessel functions of the first kind $(J_n)$ times a Hankel function of the first kind $(H_n^{(1)})$ when there are constraints between the three orders and the power-law weight that enable factorizing the Appell F4 function of the free arguments (appearing in Bailey's formula, \cite{bailey}) into a product of two hypergeometrics; interestingly, this work finds a triangular inequality between the three free arguments, as often appears in these types of integrals. \cite{gervois86} extends the approach of the previous work to modified Bessel functions of the first and second kind (respectively $I_n$ and $K_n$) and to non-triangular combinations of free arguments. \cite{whelan93} computes a three-SBF integral with $k^2$ weight and subject to the restriction that the orders satisfy a triangle inequality. He also discusses a factor of $1/2$ that must be introduced at the endpoints of the integration domain that stem from use of a Delta function; we will encounter a similar such factor in this work.\footnote{We note a typo in his result, equation 9; his $\theta(k, k', n)$ there should be $\theta(k, k', \rho)$.} 

\cite{kamion2000} uses an approach first outlined in \cite{friedberg} to recursively obtain triple-sBF integrals with a $k^2$ weight subject to the conditions that the sum of the orders is even and that the three free arguments form a triangle. The recursion they find (their equation 47) is able to raise the order of one sBF in the integrand at a time. \cite{chen04} computes three-SBF integrals with $k^2$ weight (their equation A18 and following), using contour integration and proving that the integral vanishes if the three free arguments cannot form a triangle. They also obtain a four-SBF integral with the same weight (their equation B11) and subject to the restrictions that its orders have an even sum; they show there are also several inequalities on the free arguments. Their result involves an infinite sum. \cite{shellard} presents solutions for integrals of three SBFs (see their equation 66 and onward) weighted by powers of $k^2$ and $k^{-1}$, representing the SBFs in terms of complex exponentials, giving a result in terms of a finite series (but with many terms, which aspect the authors note poses challenges to convergence when used in numerical work). Their procedure is a development of one originally presented in \cite{gervois89}, which latter work they note did not handle completely the analysis of singularities in the intermediate steps. There is no restriction on the sum of the orders in \cite{shellard}. 

\subsection{Most Related Works, Other Works of Interest, \& Plan of the Paper}
Perhaps the two most closely related works to the present one are, first, \cite{maximon}, which explores integrals of the form $\int k^2 dk\; j_{\ell_1}(kr_1) j_{\ell_2}(kr_2)$ where $\ell_1$ and $\ell_2$ need not be equal. That work finds that for even-parity $\ell_1 + \ell_2$, the integrals can be expressed in terms of Dirac delta functions, but for odd-parity $\ell_1 + \ell_2$, the results never involve Delta functions, but only Heavisides and kernels coupling $r_1$ and $r_2$ that are related to the Legendre functions. In the present work, we expand this result to encompass all positive-semi-definite powers of $k$ in the integrand. 

Second, \cite{stuart} uses an expansion of a spherical harmonic evaluated \textit{on} the gradient operator, as well as a translation formula, to show that double SBF integrals  can be obtained as derivatives of a Dirac delta function, subject to the condition that the sum of the power-law and the two SBF orders be even (that work's equation (26), where $\sigma$ denotes a spherically-symmetric Delta function). This is a restriction the present work will relax. Here we also treat integrals of three or more SBFs weighted by an arbitrary (integer, positive-semi-definite) power, extending Mehrem's work in that regard. Indeed, we show how the results we obtain allow one to evaluate an integral of an arbitrary number of SBFs weighted by an arbitrary power. 

Before turning to the main text, we mention briefly several other less-directly related but interesting works. \cite{dominici_id} finds an exact integral-to-sum identity for double SBF overlap integrals, and we obtain a simpler proof of this  in \cite{philcox_to_sum} as well as extending it to products of an arbitrary number of SBFs. \cite{majumdar} relates certain such integrals (\textit{e.g.} one over a product of $j_0$'s) to the behavior of random walkers, leading to interesting and simple means of proving some otherwise unintuitive results. \cite{fabrikant} adopts a parametric-differentiation approach similar to ours (but using the Rayleigh formula) to evaluate triple and quadruple SBF integrals, yet does not treat divergent cases. \cite{assassi} suggests using the spherical Bessel differential equation to raise the power of the integrand, as we do here. \cite{adkins_ft} derives a number of useful Fourier Transform pairs using single and double sBF integrals as well as delta functions and expansions into spherical harmonics. 

The paper is structured as follows. In \S\ref{sec:over}, we summarize the properties of SBFs we will find most germane. In \S\ref{sec:method} we outline our method and the base results needed. \S\ref{sec:deriv} presents our detailed calculations and final formula. We then display a few test cases in \S\ref{sec:tests}. \S\ref{sec:triple} extends the work to triple-SBF integrals, and then \S\ref{sec:more} shows how to treat an arbitrary number. We conclude in \S\ref{sec:concs}. 


\section{Brief Review of Relevant Properties}
\label{sec:over}
The spherical Bessel functions are solutions to Bessel's differential equation:
\begin{align}
    {x^2}{\frac{d^2f}{dx^2}}+{2x}{\frac{df}{dx}}+{[{x^2}-{\ell}{(\ell+1)}]f}=0,
\end{align}
where the solution $f(x) = j_{\ell}(x)$ is the spherical Bessel of the first kind, and the solution $f(x) = y_{\ell}(x)$ is the spherical Bessel of the second kind.

Spherical Bessel functions of the first kind are related to Bessel functions by
\begin{align}
\label{eqn:sbf_to_bf}
    j_{\ell}(x) = \sqrt{\frac{\pi}{2x}} J_{\ell + 1/2}(x).
\end{align}
The first few spherical Bessel functions of the first kind are:
\begin{align}
    j_{0}(x) &= {\frac{\sin{x}}{x}} =  \mathrm{sinc{(x)}}, \nonumber \\
    j_{1}(x) &= {\frac{\sin{x}}{x^2}} - {\frac{\cos{x}}{x}}, \nonumber \\
    j_{2}(x) &= \left(\frac{3}{x^2}-1\right)\frac{\sin{x}}{x}-\frac{3\cos{x}}{x^2}.
\end{align}
Rayleigh's formula generates the order-$\ell$ SBF from the order-zero one, as
\begin{align}
     j_{\ell}(x) = {(-x)^{\ell}}{\left(\frac{1}{x}\frac{d}{dx}\right)}^{\ell}{\frac{\sin{x}}{x}}.
\end{align}
The behavior for small values of $x$ is
\begin{align}
    j_{\ell}(x) \approx {\frac{x^{\ell}}{(2{\ell}+1)!!}},
\end{align}
where $!!$ denotes the double factorial, which is the product of all integers from one up to $n$ (inclusive) that share the parity of $n$. The behavior for large values of $x$ is
\begin{align}
\label{eqn:large_asymp}
    j_{\ell}(x) \approx {\frac{1}{x}}{\sin{\left(x-{\frac{{\ell}{\pi}}{2}}\right)}}.
\end{align}
The closure relation is:
\begin{align}
\label{eqn:close}
    \int_{0}^{\infty}x^2 dx\; j_{\ell}(kx) j_{\ell}(k'x) = {\frac{\pi}{2kk'}}{ \delta_{\rm D}^{[1]}}(k-k'),
\end{align}
where $\delta_{\rm D}^{[1]}$ is a 1D Dirac delta function.
\\
\\The orthogonality relation is:
\begin{align}
    \int_{-\infty}^{\infty}{j_{\ell}(x)}{j_{\ell'}(x)}{dx} = {\frac{\pi}{2{\ell}+1}}{\delta^{\rm K}_{\ell \ell'}},
\end{align}
where $\delta^{\rm K}_{\ell \ell'}$ is the Kronecker delta, which is unity when $\ell = \ell'$ and zero otherwise.

\section{Method and Fundamental Results Used}
\label{sec:method}
\subsection{Integral Form and Operator}
We are concerned with integrals of the form:
\begin{align}
\label{eqn:GenInt}
    I_{\ell \ell'}^{[n]}(r, r') = \int_0^{\infty} dk k^n j_{\ell}(kr) j_{\ell'}(kr').
\end{align}
In what follows, we suppress the upper and lower bounds as we will always be considering half-infinite integrals.

We may use Bessel's differential equation to define an operator that will act on the SBF and by parametric differentiation introduce an additional power of $k^2$:
\begin{align}
\label{eqn:OpEq}
    \hat{D}_{\ell, \chi} j_{\ell}(k\chi) \equiv -\left[\frac{\partial^2}{\partial \chi^2} + \frac{2}{\chi} \frac{\partial}{\partial \chi}  - \frac{\ell(\ell+1)}{\chi^2}\right]\;j_{\ell}(k\chi) = k^2 j_{\ell}(k\chi). 
\end{align}
This method was suggested and used to write a formal result for even powers of $k$ against two SBFs of the same order in \cite{assassi}.

To use this operator, we need a base integral for equation (\ref{eqn:GenInt}) for some value of $n$. There are certain $n$ for which the integral is convergent; we begin with them. We will need a formula that has both an odd and an even allowed starting $n$ since applying our operator cannot alter the parity of the power of $k$.

We observe that
\begin{align}
     \hat{D}_{\ell, r}I_{\ell \ell'}^{[n]}(r, r') = I_{\ell \ell'}^{[n+2]}(r, r') = \int dk\, k^{n+2} j_{\ell}(kr) j_{\ell'}(kr').
\end{align}

The $\hat{D}_{\ell,r}$ operator can be repeatedly applied to form a ``ladder'' of integrals with $n$ ascending in steps of two. This means that an integral of the form (\ref{eqn:GenInt}) can be evaluated for arbitrary $n$ by applying the $\hat{D}$ operator a suitable number of times, as long as we have a result for some starting value of $n$; in particular we will need both even and odd parity-$n$ starting cases as the operator cannot change the parity of the power.

\subsection{Base Case}
We find in \cite{GR07} integrals that offer the required base cases for even and odd $n$. We set variables in their work to those in ours as follows: $\alpha \to r$, $\beta \to r'$, $\nu \to \ell+1/2$, $\mu \to \ell'+1/2$, $t \to k$, $\lambda \to 1-n$. We use colors to separate the cases in \cite{GR07}; \textcolor{blue}{blue} is for \textcolor{blue}{$r'>r$} and \textcolor{red}{red} is for \textcolor{red}{$r>r'$}. We now present the base results we will require, showing both the full form and then a shorthand one in terms of coefficients which in each case we will use in what follows. We use equation (\ref{eqn:sbf_to_bf}) to convert the spherical Bessel functions into Bessel functions; the \cite{GR07} formulae we use are in terms of these latter.
\begin{itemize}
        \item \textcolor{blue}{}Case where $r'>r$ (\cite{GR07} equation 6.574.1):
\end{itemize}
\begin{align}
         &\textcolor{blue}{I_{r'>}^{[n]} = \frac{\pi}{2^{2-n}}\frac{r^{\ell}}{r'^{\ell+1+n}}\frac{\Gamma\left(\frac{1}{2}(\ell+\ell'+1+n)\right)}{\Gamma\left(\frac{1}{2}(\ell'-\ell+2-n)\right)\Gamma\left(\ell+\frac{3}{2}\right)}} \nonumber \\ 
         & \qquad\textcolor{blue}{\times\;_2F_1\left(\frac{1}{2}(\ell+\ell'+n+1),\frac{1}{2}(\ell-\ell'+n);\ell+\frac{3}{2};
     \left({\frac{r}{r'}}\right)^2\right).}
     \end{align}
We define both a constant for the term in the first line that involves Gamma functions, as 
\begin{align}
\label{eqn:cgrp}
    \textcolor{blue}{c_{\Gamma r'>} = \frac{\Gamma\left(\frac{1}{2}(\ell+\ell'+1+n)\right)}{\Gamma\left(\frac{1}{2}(\ell'-\ell+2-n)\right)\Gamma\left(\ell+\frac{3}{2}\right)},}
\end{align}
and shorthand notation for the hypergeometric as:
\begin{align}
     \textcolor{blue}{F_{'xyz} =\; _2F_1\left(\frac{1}{2}(\ell+\ell'+n+x),\frac{1}{2}(\ell-\ell'+n+y);\ell+\frac{z}{2};
     \left({\frac{r}{r'}}\right)^2\right)}.
     \label{eqn:sh_blue}
\end{align}

Using our definitions to shorten the coefficients and the hypergeometric, this becomes
     \begin{align}
    &\textcolor{blue}{\boxed{I_{r'>}^{[n]} = \frac{\pi}{2^{2-n}}c_{\Gamma r'>}\frac{r^\ell}{r'^{\ell+1+n}}F_{'103}}}
    \label{eqn:rp_big}
\end{align}
where $\ell+\ell'+n+1>0$, $n<2$, and $0<r<r'$.

\begin{itemize}
        \item \textcolor{red}{}Case where $r' < r$ (\cite{GR07} equation 6.574.3)
\end{itemize}
\begin{align}
    &\textcolor{red}{I_{r>}^{[n]} = \frac{\pi}{2^{2-n}}\frac{r'^{\ell}}{r^{\ell'+1+n}}\frac{\Gamma\left(\frac{1}{2}(\ell+\ell'+1+n)\right)}{\Gamma\left(\frac{1}{2}(\ell-\ell'+2-n)\right)\Gamma\left(\ell+\frac{3}{2}\right)}} \nonumber \\ &\qquad\textcolor{red}{\times\;_2F_1\left(\frac{1}{2}(\ell+\ell'+n+1),\frac{1}{2}(\ell'-\ell+n);\ell'+\frac{3}{2};
     \left({\frac{r'}{r}}\right)^2\right).} \nonumber
     \end{align}
We again define a constant for the term in the first line involving the Gamma functions:
\begin{align}
\label{eqn:cgr}
    \textcolor{red}{c_{\Gamma r>} = \frac{\Gamma\left(\frac{1}{2}(\ell+\ell'+1+n)\right)}{\Gamma\left(\frac{1}{2}(\ell-\ell'+2-n)\right)\Gamma\left(\ell+\frac{3}{2}\right)},}
\end{align}
and a shorthand for the hypergeometric:
\begin{align}
    \textcolor{red}{F_{xyz} = \;     _2F_1\left(\frac{1}{2}(\ell+\ell'+n+x),\frac{1}{2}(\ell'-\ell+n+y);\ell'+\frac{z}{2};
     \left({\frac{r'}{r}}\right)^2\right)}.
     \label{eqn:sh_red}
\end{align}
Using our definitions to shorten the coefficients and the hypergeometric, this becomes
     \begin{align}
    \textcolor{red}{\boxed{I_{r>}^{[n]}=\frac{\pi}{2^{2-n}}c_{\Gamma r>}\frac{r'^\ell}{r^{\ell'+1+n}}F_{103}}}
    \label{eqn:rbig}
\end{align}
where $\ell+\ell'+n+1>0$, $n<2$, and $0<r'<r$.

Instead of applying the operator to each of these cases separately, we combine them into one expression for the integral to which the operator may then be applied. We do so using Heaviside functions, as
\begin{align}
\label{eqn:IntHK}
    I_{\ell \ell'}^{[n]}(r, r') =\textcolor{blue}{H(r'-r)I_{r'>}}+\textcolor{red}{H(r-r')I_{r>}}.
\end{align}
We intend the Heaviside functions as left-continuous $(H(0)=0)$ here. In our representation (\ref{eqn:IntHK}), we omit the value of $I_{\ell\ell'}^{[n]}$ at $r=r'$, which is given by \cite{GR07} equation 6.574.2, and is not a limit of the other two cases (as when $r \to r'$ the argument of the hypergeometric for them becomes unity, which can cause the hypergeometric to diverge). We omit this value because it is only relevant on a set of measure zero. Hence it is not needed if the results obtained in this work are integrated over, as we expect would be the case in applications. Indeed, the use of Dirac delta functions is justified only if they are ultimately integrated over. In contrast to the $r=r'$ value, we note that Delta functions, defined for instance as the limit of a Gaussian as its width goes to zero, are always supported on a set of nonzero measure even as this limit is taken; hence they must be retained in the calculation. 

\section{Double-SBF Overlap Integral General Formula}
\label{sec:deriv}
The $\hat{D}$ operator outlined in equation (\ref{eqn:OpEq}) can now be applied to the RHS of equation (\ref{eqn:IntHK}). We preserve our organization by colors by maintaining a term's color when a derivative is applied.

Beginning with the second term of the operator in equation (\ref{eqn:OpEq}) and for the moment focusing only on the derivative and leaving the factor of $2/r$ unapplied (as we will use this raw first derivative to obtain the second derivative shortly), we differentiate equation (\ref{eqn:IntHK}). We use the product rule, the chain rule, and the relation $dH(x)/dx = \dD(x)$, where $\dD(x)$ is the 1D Dirac delta function:
\begin{align}
\label{eqn:firstder}
    \frac{\partial}{\partial r}[I_{\ell \ell'}^{[n]}(r, r')] = \textcolor{blue}{-\frac{1}{2}\dD(r'-r)I_{r'>}+H(r'-r)\frac{\partial}{\partial r}I_{r'>}} \nonumber \\ +\textcolor{red}{\frac{1}{2}\dD(r-r')I_{r>}+ H(r-r')\frac{\partial}{\partial r}I_{r>}}.
\end{align}
We retain the color of signs when they originate from applying the chain rule.

We now pause to outline an important subtlety: a careful reader will have noted the pre-factors of $1/2$ introduced by hand on the Dirac Delta functions. This is because in fact we wish to take only left-handed or right-handed derivatives of the Heaviside functions; the \textcolor{blue}{blue} region derivative may only be taken as $r$ approaches $r'$ from the left, whereas the \textcolor{red}{red} region derivative may only be taken as $r$ approaches $r'$ from the right. The symmetry of the Dirac delta function implies that each of these ``handed'' derivatives must be half of the total derivative. This argument is along the same lines as that in \cite{whelan93} where he argues that a Dirac delta function evaluated through the endpoint of an interval should only be counted with $1/2$ weight (comment below his equation 8). 

We now differentiate equation (\ref{eqn:firstder}) to yield the first term of the operator in equation (\ref{eqn:IntHK}):
\begin{align}
\label{eqn:secondder}
    \frac{\partial^2}{\partial r^2}I_{\ell \ell'}^{[n]}(r, r') &=     \textcolor{blue}{-\frac{\partial}{\partial r}\left[\frac{1}{2}\dD(r'-r)\right]I_{r'>} - \frac{1}{2}\dD(r'-r)\frac{\partial}{\partial r}I_{r'>}} \nonumber \\
    &\textcolor{blue}{- \frac{1}{2}\dD(r'-r)\frac{\partial}{\partial r}I_{r'>} + H(r'-r)\frac{\partial^2}{\partial r^2}I_{r'>}} \nonumber \\ 
    &+\textcolor{red}{\frac{\partial}{\partial r}\left[\frac{1}{2}\dD(r-r')\right]I_{r>} + \frac{1}{2}\dD(r-r')\frac{\partial}{\partial r}I_{r>}} \nonumber \\ &\textcolor{red}{+\frac{1}{2}\dD(r-r')\frac{\partial}{\partial r}I_{r>} + H(r-r')\frac{\partial^2}{\partial r^2}I_{r>}}. 
    \end{align}
Now substituting equations (\ref{eqn:firstder}) and (\ref{eqn:secondder}) into equation (\ref{eqn:IntHK}) for the operator, we find:
\begin{align}
    &I_{\ell \ell'}^{[n+2]}(r, r') = -\left[\frac{\partial^2}{\partial r^2} + \frac{2}{r} \frac{\partial}{\partial r}  - \frac{\ell(\ell+1)}{r^2}\right]I_{\ell \ell'}^{[n]}(r, r') \nonumber \\
     & \qquad = - \bigg\{ \bigg[\textcolor{blue}{-\frac{\partial}{\partial r}\left[\frac{1}{2}\dD(r'-r)\right] I_{r'>} - \dD(r'-r)\frac{\partial}{\partial r}I_{r'>} + H(r'-r)\frac{\partial^2}{\partial r^2}I_{r'>}} \nonumber \\ 
    & \qquad + \textcolor{red}{\frac{\partial}{\partial r}\left[ \frac{1}{2} \dD(r-r')\right] I_{r>} + \dD(r-r')\frac{\partial}{\partial r}I_{r>}+ H(r-r')\frac{\partial^2}{\partial r^2}I_{r>}}\bigg] \nonumber \\
    & \qquad +\frac{2}{r} \bigg[    - \textcolor{blue}{\frac{1}{2}\dD(r'-r)I_{r'>}
    +H(r'-r)\frac{\partial}{\partial r}I_{r'>}} + \textcolor{red}{\frac{1}{2} \dD(r-r')I_{r>}+ H(r-r')\frac{\partial}{\partial r}I_{r>}}\bigg] \nonumber \\
    &\qquad -\bigg[ \frac{\ell(\ell +1)}{r^2} \Big(\textcolor{blue}{H(r'-r)I_{r'>}}+\textcolor{red}{H(r-r')I_{r>}}
    \Big)\bigg]\bigg\}.
\end{align}
We now group terms by the region in $(r,r')$ space on which they will apply, while retaining their original color. We introduce a \textcolor{ForestGreen}{green} bracket to group terms that will apply when \textcolor{ForestGreen}{$r=r'$}. We also factor out the Delta functions and their derivatives by using their even parity, \textit{i.e.} that $\dD(r-r')=\dD(r'-r)$, and also the relation ${\partial}\,\dD(r-r')/{\partial r} = {\partial}\,\dD(r'-r)/{\partial r}$. We obtain:
\begin{align}
\label{eqn:ExprsansIs}
    &I_{\ell \ell'}^{[n+2]}(r, r') = -\bigg\{\textcolor{blue}{H(r'-r)\bigg(\frac{\partial^2}{\partial r^2}I_{r'>}+\frac{2}{r}\frac{\partial}{\partial r}I_{r'>}-\frac{\ell(\ell+1)}{r^2}I_{r'>}\bigg)} \nonumber \\ &+\textcolor{red}{H(r-r')\bigg(\frac{\partial^2}{\partial r^2}I_{r>}+\frac{2}{r}\frac{\partial}{\partial r}I_{r>}-\frac{\ell(\ell+1)}{r^2}I_{r>}\bigg)}\nonumber \\
    &+ \textcolor{ForestGreen}{\bigg[}
    \dD(r-r')\bigg(\textcolor{blue}{-\frac{\partial}{\partial r}I_{r'>}-\frac{I_{r'>}}{r}}+
    \textcolor{red}{\frac{\partial}{\partial r}I_{r>}+\frac{I_{r>}}{r}}\bigg)+\frac{\partial}{\partial r} 
    \left[\frac{1}{2} \dD(r-r')\right] \bigg(\textcolor{blue}{-I_{r'>}}+\textcolor{red}{I_{r>}}\bigg) \textcolor{ForestGreen}{\bigg]}\bigg\}.
\end{align}

We now must substitute our expressions for $\textcolor{blue}{I_{r'>}}$ and $\textcolor{red}{I_{r>}}$ (respectively, equations (\ref{eqn:rp_big}) and (\ref{eqn:rbig})) and their first and second derivatives into equation (\ref{eqn:ExprsansIs}). We display these latter here:
\begin{align}
    &\textcolor{blue}{\frac{\partial}{\partial r}I_{r'>}=\frac{\pi}{2^{2-n}}c_{\Gamma r'>}\bigg(\frac{2r^{\ell+1}}{r'^{\ell+3+n}}d_{'103}F_{'325}+\frac{\ell r^{\ell-1}}{r'^{\ell+1+n}}F_{'103}\bigg)},\nonumber\\
    &\textcolor{blue}{\frac{\partial^2}{\partial r^2}I_{r'>}=\frac{\pi}{2^{2-n}}c_{\Gamma r'>}\bigg(\frac{4r^{\ell+2}}{r'^{\ell+5+n}}d_{'103}d_{'325}F_{547}+\frac{(4\ell+2)r^\ell}{r'^{\ell+3+n}}d_{'103}F_{'325} +\frac{(\ell^2-\ell)r^{\ell-2}}{r'^{\ell+1+n}}F_{'103}\bigg)},\nonumber\\
    &\textcolor{red}{\frac{\partial}{\partial r}I_{r>}=-\frac{\pi}{2^{2-n}}c_{\Gamma r>}\bigg(\frac{(\ell'+1+n)r'^\ell}{r^{\ell'+2+n}}F_{103}+\frac{2r'^{\ell+2}}{r^{\ell'+4+n}}d_{103}F_{325}\bigg)},\nonumber\\
    &\textcolor{red}{\frac{\partial^2}{\partial r^2}I_{r>}=\frac{\pi}{2^{2-n}}c_{\Gamma r>}\bigg(\frac{(\ell'+2+n)(\ell'+1+n)r'^\ell}{r^{\ell'+3+n}}F_{103}+\frac{(4\ell'+10+4n)r'^{\ell+2}}{r^{\ell'+5+n}}d_{103}F_{325}}\nonumber \\
    &\qquad \qquad\;\textcolor{red}{+\frac{4r'^{\ell+4}}{r^{\ell'+7+n}}d_{103}d_{325}F_{'547}\bigg)}.
\end{align}
We have used the following relation for the derivative of a hypergeometric function:\footnote{\url{https://functions.wolfram.com/HypergeometricFunctions/Hypergeometric2F1/20/01/05/}, first formula.}
\begin{align}
    \frac{\partial}{\partial u} \; _2F_1(a,b;c; u) = \frac{ab}{c} \; _2F_1(a + 1, b + 1; c + 1; u).
\end{align}
Given the factors of $1/2$ on each parameter in our shorthand $F_{xyz}$ (equations \ref{eqn:sh_blue} and \ref{eqn:sh_red}), raising each parameter by unity raises each lower index by two, so the hypergeometric factor in the first derivative $F_{103}$ becomes $F_{325}$.

We have also defined a shorthand for the constant $(ab)/c$ that stems from differentiating the hypergeometric function:
\begin{align}
\label{eqn:d'xyz}
    \textcolor{blue}{d_{'xyz} = \frac{(\ell+\ell'+n+x)(\ell-\ell'+n+y)}{4(\ell+{z}/{2})}},
\end{align}
and
\begin{align}
\label{eqn:dxyz}
    \textcolor{red}{d_{xyz} = \frac{(\ell+\ell'+n+x)(\ell'-\ell+n+y)}{4(\ell'+{z}/{2})}};
\end{align}
the $4$ in the denominators comes from the fact that the first two parameters of the hypergeometric have overall factors of $1/2$ that may be multiplied together.

Finally, substituting these expressions into equation (\ref{eqn:ExprsansIs}), simplifying and factoring out constants, we have:
\begin{align}
    &I_{\ell \ell'}^{[n+2]}(r, r') =  -\frac{\pi}{2^{2-n}} \Bigg\{ \textcolor{blue}{H(r'-r)c_{\Gamma r'>} \bigg(\frac{4r^{\ell+2}}{r'^{\ell+5+n}}d_{'103}d_{'325}F_{'547}+\frac{(4\ell+2)r^\ell}{r'^{\ell+3+n}}d_{'103}F_{'325}} \nonumber \\ &\;\;\;\textcolor{blue}{+\frac{(\ell^2-\ell)r^{\ell-2}}{r'^{\ell+1+n}}F_{'103}+\left(\frac{4r^\ell}{r'^{\ell+3+n}}d_{'103}F_{'325}+\frac{2\ell r^{\ell-2}}{r'^{\ell+1+n}}F_{'103}\bigg)-\frac{\ell(\ell+1)r^{\ell-2}}{r'^{\ell+1+n}}F_{'103}\right)} \nonumber \\
    &\;\;\;+ \textcolor{red}{ H(r-r')c_{\Gamma r>}\bigg(\frac{(\ell'+2+n)(\ell'+1+n)r'^\ell}{r^{\ell'+3+n}}F_{103}+\frac{(4\ell'+10+4n)r'^{\ell+2}}{r^{\ell'+5+n}}d_{103}F_{325}}\nonumber \\
    &\;\;\;\textcolor{red}{+\frac{4r'^{\ell+4}}{r^{\ell'+7+n}}d_{103}d_{325}F_{547}+\left(-\frac{2(\ell'+1+n)r'^\ell}{r^{\ell'+3+n}}F_{103}-\frac{4r'^{\ell+2}}{r^{\ell'+5+n}}d_{103}F_{325}\right)} 
    \textcolor{red}{-\frac{\ell(\ell+1)r'^\ell}{r^{\ell+3+n}}F_{103}\bigg)} \nonumber \\
    &\;\;\;+ \textcolor{ForestGreen}{\Bigg[}\dD(r-r')\bigg(\textcolor{blue}{c_{\Gamma r'>}\left(-\frac{\ell r^{\ell-1}}{r'^{\ell+1+n}}F_{103}-\frac{2r^{\ell+1}}{r'^{\ell+3+n}}d_{'103}F_{'325}-\frac{r^{\ell-1}}{r'^{\ell+1+n}}F_{'103}\right)} \nonumber \\
    &\;+ \textcolor{red}{c_{\Gamma r>}\left(-\frac{(\ell'+1+n)r'^\ell}{r^{\ell'+2+n}}F_{103}-\frac{2r'^{\ell+2}}{r^{\ell'+4+n}}d_{103}F_{325}+\frac{r'^\ell}{r^{\ell'+2+n}}F_{103}\right)}\bigg) \nonumber \\
    &\;\;\;+\frac{\partial}{\partial r}\left[ \frac{1}{2} \dD(r-r')\right] \left(\textcolor{blue}{-\frac{c_{\Gamma r>}r^\ell}{r'^{\ell+1+n}}F_{'103}}+\textcolor{red}{\frac{c_{\Gamma r'>}r'^\ell}{r^{\ell'+1+n}}F_{103}}\right)
    \textcolor{ForestGreen}{\bigg]}\Bigg\}. 
\end{align}
Combining terms we find:
\begin{align}
\label{eqn:FinalExpression}
     &I_{\ell \ell'}^{[n+2]}(r, r') =  -\frac{\pi}{2^{2-n}} \Bigg\{ \textcolor{blue}{H(r'-r)c_{\Gamma r'>} \bigg(\frac{4r^{\ell+2}}{r'^{\ell+5+n}}d_{'103}d_{'325}F_{'547}+\frac{(4\ell+6)r^\ell}{r'^{\ell+3+n}}d_{'103}F_{'325}\bigg)} \nonumber \\
     &\;+ \textcolor{red}{ H(r-r')c_{\Gamma r>}\bigg(\frac{\left(2n\ell'+\ell'+n^2+\ell'^2-\ell^2-\ell + n\right)r'^\ell}{r^{\ell'+3+n}}F_{103}+\frac{(4\ell'+6+4n)r'^{\ell+2}}{r^{\ell'+5+n}}d_{103}F_{325}} \nonumber \\
     &\;\textcolor{red}{+\frac{4r'^{\ell+4}}{r^{\ell'+7+n}}d_{103}d_{325}F_{547}\bigg)} \nonumber \\
     &\; + \textcolor{ForestGreen}{\Bigg[}\dD(r-r')\bigg(\textcolor{blue}{c_{\Gamma r'>}\left(-\frac{2r^{\ell+1}}{r'^{\ell+3+n}}d_{'103}F_{'325}-\frac{(\ell+1) r^{\ell-1}}{r'^{\ell+1+n}}F_{'103}\right)} \nonumber \\
    &\; + \textcolor{red}{c_{\Gamma r>}\left(-\frac{2r'^{\ell+2}}{r^{\ell'+4+n}}d_{103}F_{325}-\frac{(\ell'+n)r'^\ell}{r^{\ell'+2+n}}F_{103}\right)}\bigg) \nonumber \\
    &\; +\frac{\partial}{\partial r}\left[ \frac{1}{2} \dD(r-r') \right]{}\left(\textcolor{blue}{-\frac{c_{\Gamma r'>}r^\ell}{r'^{\ell+1+n}}F_{'103}}+\textcolor{red}{\frac{c_{\Gamma r>}r'^\ell}{r^{\ell'+1+n}}F_{103}}\right)    
     \textcolor{ForestGreen}{\Bigg]}\Bigg\}.
\end{align}
This is our final expression for the integral form outlined in equation (\ref{eqn:GenInt}), and is a key result of this work. We have used \textsc{Mathematica} to confirm the reduction starting directly from equation (\ref{eqn:IntHK}).\footnote{Relative to the work in \textsc{Mathematica}, we note three points. First, when both $r$ and $r'$ are written purely symbolically (\textit{i.e.} without any values inserted for either, especially $r'$), \textsc{Mathematica} writes the Heaviside function's second derivative simply as ${\rm DiracDelta'[r-r']}$ or ${\rm DiracDelta'[r'-r]}$, where prime is ambiguous and this result does not capture the negative sign produced on differentiating the Heaviside $H(r' - r)$ for the first time (by chain rule). In contrast, \textsc{Mathematica} does perform the differentiation correctly if a numerical value is specified for $r'$. \textsc{Mathematica} also performs a first differentiation of the Heaviside correctly even with both variables' remaining symbolic; the issue only appears when the second derivative is requested \textit{ab initio}. Thus, the sign of the first blue term in the last line of equation (\ref{eqn:FinalExpression}) requires some thought to compare with the \textsc{Mathematica} result; the negative sign we have above is intended. Second, the factors of $1/2$ on the Delta functions and their derivatives (as discussed around equation (\ref{eqn:firstder})) are not present in the \textsc{mathematica} notebook as it does not handle correctly the question of left and right-handed derivatives. Third, \textsc{mathematica} often uses the regularized hypergeometric$\;_2 \tilde{{\rm F}}_1$; a helpful relation is that$ \;_2 \tilde{{\rm F}}_1(a,b,c;u) = [1/\Gamma(c)]\;_2{\rm F}_1(a,b,c;u)$.}  


\section{Tests}
\label{sec:tests}
Let us begin by observing the behavior of some of the constants that we have defined. The constants $d_{xyz}$ and $d_{'xyz}$ that arose from differentiating the hypergeometric function defined in equations ({\ref{eqn:d'xyz}}) and (\ref{eqn:dxyz})  will both become zero when $n+y=0$ and $\ell=\ell'$ for any $\ell$ and $\ell'$. Additionally, $d_{xyz}=d_{'xyz}$ when $\ell=\ell'$. Finally, we note that when $\ell=\ell'$, $c_{\Gamma r>}=c_{\Gamma r'>}$.

\subsection{$\ell = 0 = \ell',\; n = 0$}
\label{subsec:delta_test}
This test will ultimately yield an integral $\int k^2 dk\;j_0(kr) j_0(kr')$, which we know gives a spherically-symmetric Dirac delta function by the closure relation (\ref{eqn:close}).

We will need to substitute our test parameters into equation (\ref{eqn:FinalExpression}). We first observe from equations (\ref{eqn:dxyz}) and (\ref{eqn:d'xyz}) that $d_{103} = 0  = d_{'103}$ in this case, meaning all terms in equation (\ref{eqn:FinalExpression}) containing either of these coefficients will immediately drop out. Any terms proportional solely to sums of powers of $\ell$, $\ell'$, and $n$, will also drop out. We find
\begin{align}
    &I_{0 0}^{[0+2]}(r, r') =  -\frac{\pi}{4} \Bigg\{ \textcolor{ForestGreen}{\Bigg[}\dD(r-r')\bigg(\textcolor{blue}{c_{\Gamma r'>}\left(-\frac{1}{r r'}F_{'103}\right)}\nonumber\\
    &\; +\frac{\partial}{\partial r}\left[ \frac{1}{2} \dD(r-r')\right]\left(\textcolor{blue}{-\frac{c_{\Gamma r'>}}{r'}F_{'103}}+\textcolor{red}{\frac{c_{\Gamma r>}}{r}F_{103}}\right)\textcolor{ForestGreen}{\Bigg]}\Bigg\}.  
\end{align}
We see that for this case, from equations (\ref{eqn:cgr}) and (\ref{eqn:cgrp}) we have $c_{\Gamma r>}= 2 = c_{\Gamma r'>}$. Using this, we find
\begin{align}
    &I_{0 0}^{[0+2]}(r, r') =  -\frac{\pi}{4} \Bigg\{
     \textcolor{ForestGreen}{\Bigg[}\dD(r-r')\bigg(\textcolor{blue}{\left(-\frac{2}{r r'}F_{'103}\right)}\bigg) \nonumber \\
    &\; +\frac{\partial}{\partial r}\left[ \frac{1}{2} \dD(r-r')\right] \left(\textcolor{blue}{-\frac{2}{r'}F_{'103}}+\textcolor{red}{\frac{2}{r}F_{103}}\right)\textcolor{ForestGreen}{\Bigg]}\Bigg\}.  
\end{align}
Using that $\partial\;\dD(r - r')/\partial r = -(r - r') \dD(r - r')$,\footnote{By applying \url{https://mathworld.wolfram.com/DeltaFunction.html} equation 14 and chain rule.} we see that in both lines, we must have $r = r'$. In view of equations (\ref{eqn:sh_blue}) and (\ref{eqn:sh_red}), this enables us to evaluate $F_{103} = 1 = F_{'103}$, which in turn renders the second line zero. Simplifying, we find
\begin{align}
    &I_{0 0}^{[2]}(r, r') =\frac{\pi}{2 r r'}  
     \dD(r-r').
\end{align}
This recovers the closure relation (\ref{eqn:close}) as desired.

\subsection{$\ell = \ell',\; n = 0$}
Just as with the above case, since we have $\ell=\ell'$ and $n=0$, all terms in equation (\ref{eqn:FinalExpression}) containing $d_{103}$ or $d_{'103}$ will vanish. We also set $n = 0$. We find:
\begin{align}
    &I_{0 0}^{[0+2]}(r, r') =  -\frac{\pi}{4} \Bigg\{
     \textcolor{red}{H(r-r')c_{\Gamma r'>}\bigg(\frac{\left(\ell'+\ell'^2-\ell^2-\ell\right)r'^\ell}{r^{\ell'+3}}F_{103}\bigg)} \nonumber \\
     &\; + \textcolor{ForestGreen}{\Bigg[}\dD(r-r')\bigg(\textcolor{blue}{c_{\Gamma r'>}\left(-\frac{(\ell+1) r^{\ell-1}}{r'^{\ell+1}}F_{'103}\right)}
     + \textcolor{red}{c_{\Gamma r>}\left(-\frac{(\ell')r'^\ell}{r^{\ell'+2}}F_{103}\right)}\bigg) \nonumber \\
    &\; +\frac{\partial}{\partial r}\left[ \frac{1}{2} \dD(r-r')\right] \left(\textcolor{blue}{-\frac{c_{\Gamma r'>}r^\ell}{r'^{\ell+1}}F_{'103}}+\textcolor{red}{\frac{c_{\Gamma r>}r'^\ell}{r^{\ell'+1}}F_{103}}\right)\textcolor{ForestGreen}{\Bigg]}\Bigg\}.  
\end{align}
Setting $\ell = \ell'$ in the first line makes it vanish. In the second line, we use that the Delta function requires $r = r'$ and hence $F_{'103} = F_{103} = 1$, and $c_{\Gamma r'>} = c_{\Gamma r>}$. In the third line, we again use $c_{\Gamma r>}  = c_{\Gamma r'>}$ and also the comment about the Delta function's derivative made in \S\ref{subsec:delta_test}, which means $r = r'$ and $F_{103} = F_{'103}$. With these substitutions the third line drops out. We have
\begin{align}
    &I_{0 0}^{[0+2]}(r, r') =  \frac{\pi}{4}  \dD(r-r')\bigg( \frac{(2\ell + 1) c_{\Gamma r>}}{r^2}\bigg).  
\end{align}
Using now that $c_{\Gamma r>}=c_{\Gamma r'>}= \Gamma(\ell+1/2)/\Gamma(\ell+3/2) = 2/(2\ell + 1)$, and replacing one $r$ with $r'$ for symmetry, we find
\begin{align}
    &I_{0 0}^{[0+2]}(r, r')  =
    \frac{\pi}{2 r r'} \dD(r-r')
\end{align}
as desired.

\subsection{$\ell = 0, \;\ell'=1,\; n = 0$}
\begin{figure}
    \centering
    \includegraphics[width=0.4\textwidth]{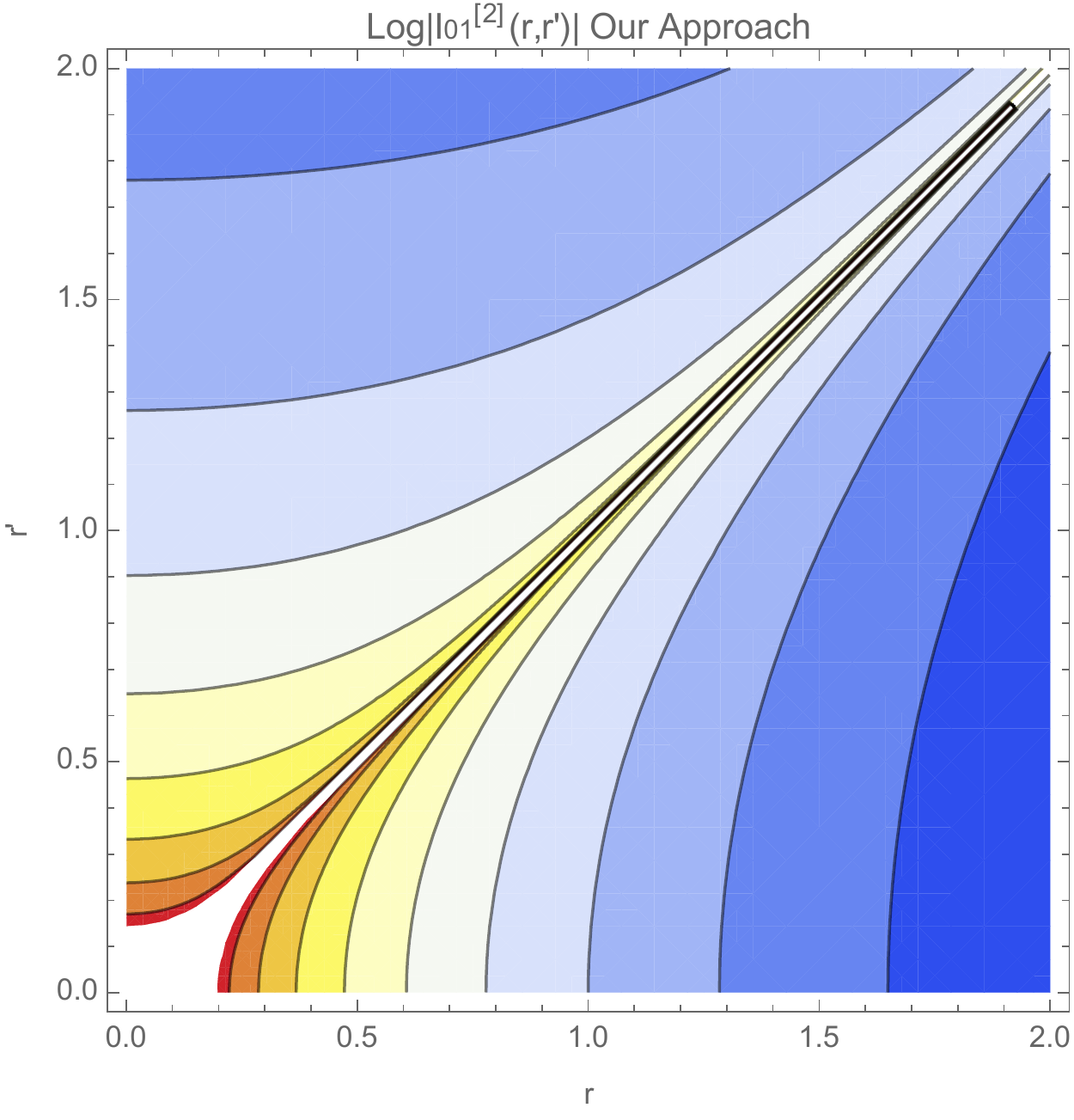}
    \includegraphics[width=0.4\textwidth]{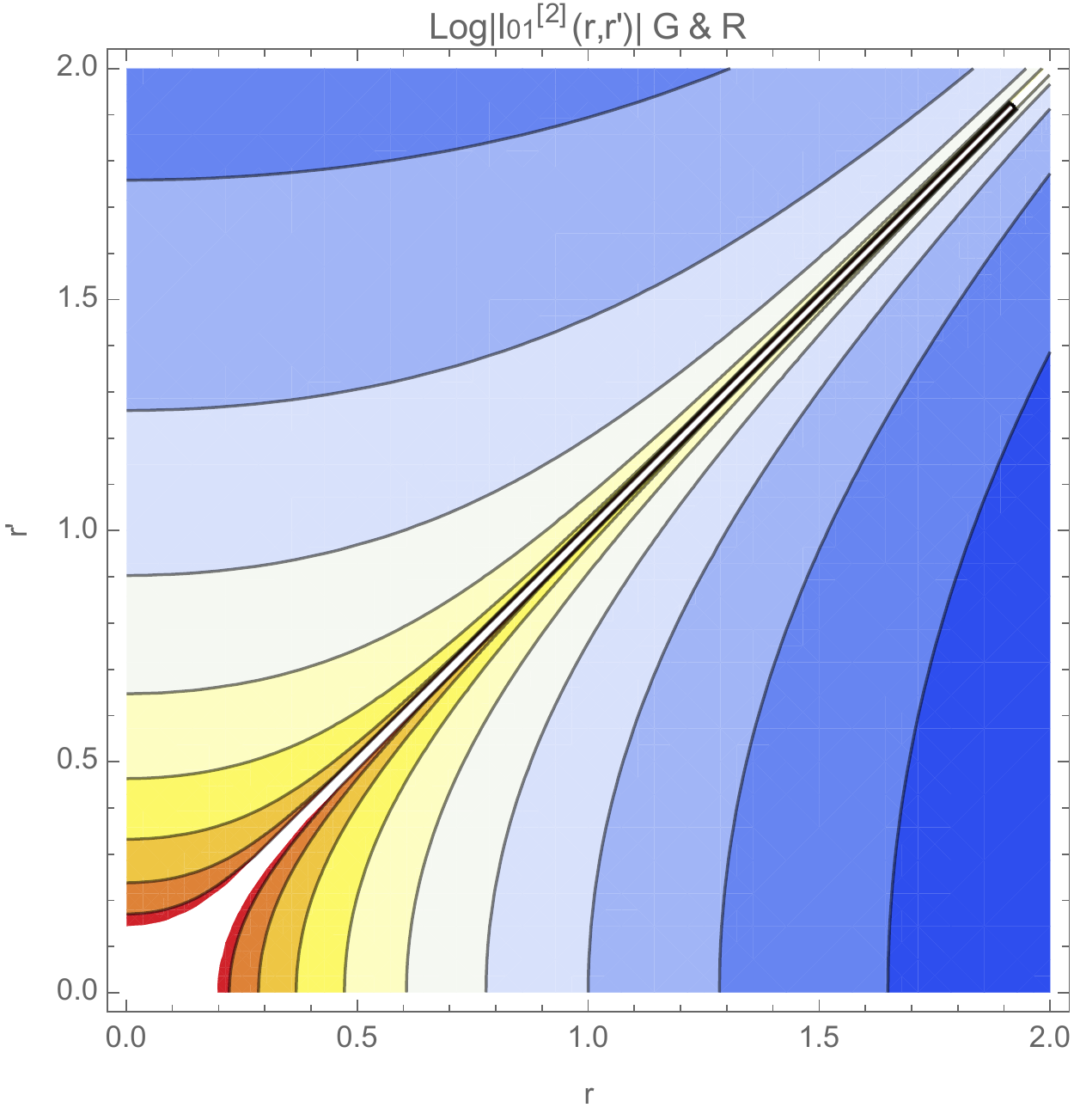}
    \includegraphics[width=0.044\textwidth]{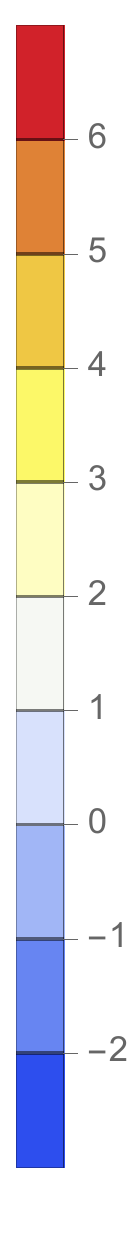}
    \caption{Comparison of using our formula \textit{(left)} vs. using Gradshteyn \& Ryzhik (\cite{GR07}) directly \textit{(right)}. The axes are $r$ (horizontal) and $r'$ (vertical) and we have taken the logarithm of the absolute value to show the dynamic range of these integrals. Overall, one can see that our approach agrees quite well with the known Gradshteyn \& Ryzhik \citep{GR07} result (given by our equation (\ref{eqn:IntHK})); we do not show a ratio to avoid dividing by numbers near zero. It is worth observing the lack of reflection symmetry about the diagonal; this is because there is no interchange symmetry between $r$ and $r'$, as we expect since they are the free arguments of sBFs of different orders (respectively, zero and one).}
    \label{fig:01_comparison}
\end{figure}

For this case, the integral resulting after raising the power is captured by equation (\ref{eqn:IntHK}) (a splicing together of two formulae from \cite{GR07}), as the result is non-singular. However, the reduction is involved, and rather than showing it, we simply allow the accompanying \textsc{Mathematica} notebook to perform it, and plot the results, comparing to the analytic expression in \cite{GR07}. We note that none of the subtleties discussed in the footnote below equation (\ref{eqn:FinalExpression}) regarding \textsc{Mathematica}'s treatment of the Heavisides or Delta functions arise, as these subtleties only enter the singular terms and there are none here. We display results for this case in Figure \ref{fig:01_comparison}.
\section{Triple-SBF Integral Result}
\label{sec:triple}
We first outline our overall approach, then present a result from the literature that we will need, and end with several examples showing the utility of our method.
\subsection{Overall Approach}
\label{subsec:overall_3}
Consider an integral of three SBFs. We have
\begin{align}
I_{\ell_1 \ell_2 \ell_3}^{[n]}(r_1, r_2, r_3) \equiv \int k^n dk\; j_{\ell_1}(kr_1) j_{\ell_2}(kr_2) j_{\ell_3}(kr_3).
\label{eqn:three}
\end{align}
If $n = 2$ and the sum $\ell_1 + \ell_2 + \ell_3$ is even, then this integral may be done and results in a polynomial, as can be seen in examples in \citep{fonseca} and as we prove in the next subsection. We seek a fully general formula. First, let us rewrite 
\begin{align}
k^{n-2} j_{\ell_3} (k r_3) = \int  dk'\;  k'^{n-2} j_{\ell_3}(k' r_3) \delta_{\rm D}^{[1]}(k - k') = \frac{2}{\pi} \int k'^n j_{\ell_3}(k' r_3) \int u^2 du \; j_{L}(u k) j_L(u k') .
\end{align}
Inserting this result in equation (\ref{eqn:three}) we have
\begin{align}
I_{\ell_1 \ell_2 \ell_3}^{[n]}(r_1, r_2, r_3) =\frac{2}{\pi}  \int u^2 du\; \int k'^n dk'\; j_{\ell_3}(k' r_3) j_L(k'u)\; \int k^2 dk\; j_{\ell_1}(kr_1) j_{\ell_2}(kr_2) j_L(ku).
\label{eqn:3_final}
\end{align}
We now perform the inner integral over $k$. Since $L$ is free, we can choose it so that $\ell_1 + \ell_2 + L$ is even and this inner integral can be done (some examples are given in \textit{e.g.} \cite{fonseca}), and results in a sum over terms each of which just involves powers of $r_1, r_2$ and $r_3$. We present and discuss the detailed general result, obtained by \cite{mehrem_2_91}, in the following subsection. We then may use the methods outlined in earlier sections of the present work to perform the integral over $k'$ in closed form. It will yield Heaviside functions, Delta functions and their derivatives, and hypergeometrics. 

We can easily integrate Heaviside functions and Delta functions and their derivatives over $u$. Since any additional dependence on $u$ from our innermost integral will be a sum of power-laws, we can also perform that integral analytically using formulae for integrals of power-laws against hypergeometrics such as NIST DLMF 15.4.1\footnote{\url{https://dlmf.nist.gov/15.14}} (the integral is just its Mellin transform) or, if a smoothing factor is needed for convergence to then have a limit taken at the end, the formula found in \textit{e.g.} Wolfram.\footnote{\url{https://functions.wolfram.com/HypergeometricFunctions/Hypergeometric2F1/21/02/02/}} This point will become more concrete when we move, shortly, to presenting several examples.

Finally, we note that this trick of introducing a Delta function is related to the method of \cite{cohl_12}, as that work's equation (2) is another form of it though one that avoids explicit appeal to the Delta function. \cite{cohl_12} exploits this trick to prove a number of interesting and complicated integrals. \cite{cohl_15} obtains further results using this idea.

\subsection{Result for Even Sum of the Orders and Square Weight}
\label{subsec:even_square}
\cite{mehrem_2_91} equation (3.21) gives the result for an even sum of the three SBF orders ($\ell_1 + \ell_2 + L = $ even) and a weight of $k^2 dk$ in the integral, as we require to perform the innermost integral in equation (\ref{eqn:3_final}). We duplicate the result here as we will need to discuss some of its properties. We have
\begin{align}
\label{eqn:mehr_result}
I_{\ell_1 \ell_2 L,\;{\rm even}}^{[2]} &= \frac{\pi \beta(\Delta)}{4 r_1 r_2 u} i^{\ell_1 + \ell_2 - L} \sqrt{2L + 1} \left( \frac{r_1}{u}\right)^{L}
\begin{pmatrix}
\ell_1 & \ell_2 & L\\
0 & 0 & 0
\end{pmatrix}^{-1}
\sum_{\mathcal{L} = 0}^{L}{{2 L} \choose {2 \mathcal{L}}}^{1/2} \left( \frac{r_2}{r_1}\right)^{\mathcal{L}}\\
&\qquad \times \sum_{\ell} (2\ell + 1) \begin{pmatrix}
\ell_1 & L - \mathcal{L} & \ell\\
0 & 0 & 0
\end{pmatrix}
\begin{pmatrix}
\ell_2 & \mathcal{L} & \ell\\
0 & 0 & 0
\end{pmatrix}
\left\{\begin{matrix}
\ell_1 & \ell_2 & L \\
\mathcal{L} & L - \mathcal{L} & \ell \\
\end{matrix}\right\}
P_{\ell}(\Delta).\nonumber
\end{align}
We may see from the inverse of the 3-$j$ symbol in the $\ell_i$ that the result of \cite{mehrem_2_91} is only valid when their sum is even; otherwise this 3-$j$ symbol vanishes and hence its inverse is ill-defined. We note that the $2 \times 3$ matrix in curly brackets is a Wigner 6-$j$ symbol (see \textit{e.g.} NIST DLMF \S34.4\footnote{\url{https://dlmf.nist.gov/34.4}}).

Equation (3.13) of \cite{mehrem_2_91} defines $\beta$ as
\begin{align}
\label{eqn:beta}
    \beta(\Delta) \equiv \bar{\theta}(1 - \Delta) \bar{\theta}(1 + \Delta)
\end{align}
and equation (3.14) of \cite{mehrem_2_91} gives the modified step function $\bar{\theta}$ as 
\begin{align}
\label{eqn:thetabar}
    \bar{\theta}(y) = \left\{0, \;y < 0;\;\;\;1/2,\;y =0;\;\;\;1, \; y>0\right\}.
\end{align}
$\Delta$ is (\cite{mehrem_2_91} equation 3.9)
\begin{align}
    \Delta \equiv \frac{r_1^2 + r_2^2 - u^2}{2 r_1 r_2}.
\end{align}
We may see from equation (\ref{eqn:mehr_result}) that the ratios $(r_1/u)$ and $(r_2/r_2)$ will only ever give rise to terms that behave like powers of each $r_i$. We may see from the explicit representation of the Legendre polynomial $P_{\ell}$ as a finite series of powers of $\Delta$ descending in same-parity steps from $\ell$ (\textit{e.g.} Wolfram equation 32\footnote{\url{https://mathworld.wolfram.com/LegendrePolynomial.html}}) that, upon expanding these powers of $\Delta$, we will simply have a sum over products of powers of the $r_i$. Consequently, our integral $I_{\ell_1 \ell_2 L,\; {\rm even}}^{[2]}$ will always be a sum of terms each of which is a product of powers of the $r_i$, as claimed. The only additional complication is $\beta$, which involves a modified step function constraining the relationship between the $r_i$.

We might worry that we could not integrate the hypergeometrics that may result from the $k'$ integral in equation (\ref{eqn:3_final}) against the result (\ref{eqn:mehr_result}) for the $k^2 dk$ integral, which will involve $\beta$. However, we have results for \textit{indefinite} integrals of the hypergeometric against power-laws, so this integral may always be performed. In particular, we have
\begin{align}
    \int z^{\alpha - 1} \;_2 {\rm F}_1(a, b, c; z) dz = \frac{z^{\alpha}}{\alpha}\;_3{\rm F}_2(a, b, \alpha; c, \alpha + 1; z)
    \label{eqn:2F1_int}
\end{align}
from \textit{e.g.} Wolfram.\footnote{\url{https://functions.wolfram.com/HypergeometricFunctions/Hypergeometric2F1/21/01/02/01/}, second formula.} As this is an indefinite integral, there is no restriction on $\alpha$ so we can also deal with negative powers if required. For some terms in the result (\ref{eqn:mehr_result}) this may be necessary as both the ratios and the Legendre polynomial in $\Delta$ bring in negative powers of $r_1, r_2$, or $u$.

Finally, we note that \cite{fabrikant}, though it claims to be very general as to the three-SBF cases it treats, actually does not indicate how to do a triple-SBF integral where the sum of the three SBF orders and the power is even; there the cosine in \cite{fabrikant}'s equation (9) vanishes and it is not well-described how to proceed. Hence, though the \cite{fabrikant} might in principle give a simpler result than \cite{mehrem_2_91} if it could be applied, in practice it is in fact the \cite{mehrem_2_91} result that is most useful for this case.

\subsection{Examples}
To demonstrate the utility of our approach, we here give a few examples that would not be obtainable by the other methods previously in the literature.
\subsubsection{$n=2,\;\ell_1 = 0= \ell_2,\; \ell_3 = 1$}
\label{subsubsec:ex_1}
We seek to evaluate the integral
\begin{align}
I_{001}^{[2]}(r_1, r_2, r_3) = \int k^2dk\;j_0(kr_1) j_0(kr_2) j_1(kr_3).   
\end{align}
We apply equation (\ref{eqn:3_final}). We first set $L=0$ in the innermost integral of that relation, and to evaluate it then use \cite{mehrem_3_gen} Appendix B, first unnumbered result below their equation B3:
\begin{align}
\label{eqn:kint}
&\int k^2 dk\; j_0(kr_1) j_0(kr_2) j_0(ku) = \frac{1}{4} \frac{\pi \beta(\Delta)}{r_1 r_2 u},\\
&\qquad \qquad \Delta = \frac{r_1^2 + r_2^2 - u^2}{2 r_1 r_2}.\nonumber 
\end{align}
Our middle integral, over $k'$, is then 
\begin{align}
    \int k'^2 dk'\;j_1(k'r_3) j_0(k'u) &=
   H(u-r_3) \left[ \frac{1}{r_3^3 - r_3 u^2} + \frac{{\rm arctanh}[r_3/u]}{r_3^2 u} \right] \nonumber\\
    &+ H(r_3 - u) \left[ \frac{2u}{3(r_3^3 - r_3 u^2)} + \frac{2\; {\rm arctanh}[u/r_3]}{3 r_3^2} \right]
    \label{eqn:kpint}
\end{align}
We note that this integral is a case that is marginally convergent as the out-of-phase behavior of the SBFs in the UV (\textit{e.g.} see equation (\ref{eqn:large_asymp})) converts the otherwise-constant asymptotic (as $k^2 \times (1/k)^2$) into an oscillating function. Hence, though it can be evaluated with our operator method, it also can be taken directly from \cite{GR07}. Nonetheless this example is still a non-trivial case where our method of reducing a triple SBF integral is required, relative to what \cite{mehrem_2_91} would be able to evaluate.

The product of equations (\ref{eqn:kint}) and (\ref{eqn:kpint}) must then be integrated against $u^2du$, as equation (\ref{eqn:3_final}) dictates, for the final result. Importantly, the factor of $1/u$ in equation (\ref{eqn:kint}) means that, effectively, we are integrating all of equation (\ref{eqn:kpint}) against a weight of $u\,du$. We also note that due to the factor of $\beta(\Delta)$ in equation (\ref{eqn:kint}), this final integral will be over a finite rather than a half-infinite interval, but since results for the indefinite integrals for even the most complicated piece in the integrand, ${\rm arctanh}$ times a power-law, are available, this is not problematic. Indeed, since ${\rm arctanh}$ may be represented as a hypergeometric$\;_2 {\rm F}_1$, equation (\ref{eqn:2F1_int}) gives the required result.

We may solve the restriction on $\Delta$ to find the bounds of integration, denoted $u_- = r_2 - r_1$ and $u_+ = r_2 + r_1$, and overall we have
\begin{align}
    I_{001}^{[2]}(r_1, r_2, r_3) =\frac{2}{\pi} \left\{H(u-r_3) I_{u>} + H(r_3 - u) I_{u<} \right\}\bigg|_{u_-}^{u_+},
\end{align}
with
\begin{align}
    &I_{u>} = -\frac{1}{2}r_3^{-1} \ln[r_3^2 - u^2]+ \frac{1}{2r_3} r_3^{-2} \left[u\; {\rm arctanh}[r_3/u] + \frac{1}{2}r_3 \ln [r_3^2 - u^2] \right],\nonumber\\
    &I_{u<} = 2\left[-\frac{u}{3 r_3} +\frac{1}{6}\ln\left[\frac{ u+ r_3}{u-r_3}\right]\right] + \frac{1}{6} r_3^{-2} \left[ 2 u^2\; {\rm arctanh} [r_3/u] + r_3\left(2u + r_3 \ln \left[ \frac{u-r_3}{u+r_3}\right]\right) \right].
\end{align}
We have dropped the contribution of measure zero when $u=r_3$ and we also do not need to track the $1/2$ in $\bar{\theta}$ (equation \ref{eqn:thetabar}) within $\beta$ (equation \ref{eqn:beta}) because it is a set of measure zero with respect to our $u$ integration. These results may now be straightforwardly evaluated at $u_+$ and $u_-$.

\subsubsection{$n=4,\;\ell_1 = 0= \ell_2,\; \ell_3 = 2$}
We now seek to evaluate the integral
\begin{align}
I_{002}^{[4]}(r_1, r_2, r_3) = \int k^4 dk\;j_0(kr_1) j_0(kr_2) j_2(kr_3).   
\end{align}
This integral is even-parity in the integrand, so the 3-$j$ symbol of \cite{mehrem_2_91} would allow it, yet the power of $k$ is too high for that work's formula to apply; the work of \cite{fabrikant} is not applicable because the integral is even-parity and also not convergent.

First, we return to equation (\ref{eqn:3_final}) and notice that we may use the same $L$ and innermost integral as we did in \S\ref{subsubsec:ex_1}. In this case our middle integral, over $k'$, is
\begin{align}
    I_{\rm mid} \equiv \int k'^4 dk'\;j_2(k'r_3) j_0(k'u) &= \hat{D}_{0,u} \int k'^2 dk'\; j_2(k' r_3) j_0(k' u) \nonumber\\
    &= \frac{3\pi}{2r_3^3} \left[ \frac{1}{u} \delta_{\rm D}^{[1]}(r_3 - u) - \frac{1}{2} \frac{\partial }{\partial u}\delta_{\rm D}^{[1]}(r_3 - u) \right].
\end{align}
We note that the operator with respect to $u$ must be applied only to this integral, and not act on the innermost integral, the result of which contains $u$-dependence; by the same token the operator must be inside the outermost integral over $u$.
Multiplying this result with equation (\ref{eqn:kint}) for our innermost integral, we have 
\begin{align}
    I_{002}^{[4]}(r_1, r_2, r_3) &= \frac{1}{2} \int u^2 du\; \frac{\beta(\Delta)}{r_1 r_2 u} I_{\rm mid} \nonumber\\
    &= \frac{3\pi}{4 r_1 r_2 r_3^3} \int_{u_-}^{u_+} u du \; \left[ \frac{1}{u} \delta_{\rm D}^{[1]}(r_3 - u) - \frac{1}{2} \frac{\partial }{\partial u}\delta_{\rm D}^{[1]}(r_3 - u) \right].
\end{align}
The integral above is trivial to perform, and we note that it will vanish unless $r_3$ satisfies a triangular inequality with $r_1$ and $r_2$; otherwise the value of $u$ picked out by the Delta functions, $u = r_3$, will not be within the domain of  integration.

\subsubsection{$n=4,\;\ell_1 = 0= \ell_2,\; \ell_3 = 1$}
We now seek to evaluate the integral
\begin{align}
I_{001}^{[4]}(r_1, r_2, r_3) = \int k^4 dk\;j_0(kr_1) j_0(kr_2) j_1(kr_3).   
\end{align}
This integral is odd-parity in the integrand and also not convergent, and again, cannot be performed using previous techniques in the literature.

In view of equation (\ref{eqn:3_final}), we observe that we may use the same innermost integral as in the previous two subsubsections, and our middle integral is
\begin{align}
    I_{\rm mid} = \int k'^4 dk'\;j_1(k'r_3)j_0(k'u).
\end{align}
It may be computed by applying our operator to the known integral below on the righthand side, as 
\begin{align}
   I_{\rm mid} = \hat{D}_{0,u} \int k'^2 dk'\;j_1(k'r_3)j_0(k'u).
\end{align}

\section{Dealing with More than Three SBFs}
\label{sec:more}
To resolve all integrals of more than three SBFs we may use a technique suggested by \cite{mehrem_4}. Consider an integral 
\begin{align}
I_{\ell_1 \ell_2 \ell_3 \ell_4}^{[n]}(r_1, r_2, r_3, r_4) \equiv \int k^n dk\; j_{\ell_1}(kr_1) j_{\ell_2}(kr_2) j_{\ell_3}(kr_3) j_{\ell_4}(kr_4).
\label{eqn:four}
\end{align}
Using the same trick as in \S\ref{sec:triple} we see that
\begin{align}
I_{\ell_1 \ell_2 \ell_3 \ell_4}^{[n]} = \frac{2}{\pi} \int u^2 du\; \int k'^n dk'\; j_{\ell_3}(k'r_3) j_{\ell_4} (k' r_4) j_L(k' u) \int k^2 dk\; j_{\ell_1}(k r_1) j_{\ell_2} (k r_2) j_L(k u).
\label{eqn:four}
\end{align}
We make two observations. First, these integrals cannot just be done using the techniques of \cite{fonseca} because {\bf i)} one of them (that over $k'$) has a power-law other than $2$, and {\bf ii)} we are not guaranteed that both sums $\ell_1 + \ell_2 + L$ and $\ell_3 + \ell_4 + L$ will be even. One might hope that at least {\bf ii)} could be dealt with, but in fact one cannot choose an $L$ that makes both of these sums even unless the sums $\ell_1 + \ell_2$ and $\ell_3 + \ell_4$ themselves have the same parity (see also discussion in Slepian and Hou in prep.). This is not guaranteed. Thus, for both reasons {\bf i)} and {\bf ii)} the technique outlined in \S\ref{sec:triple} is required to enable performing an arbitrary integral of four SBFs.

We now investigate how to split an integral of five SBFs. We have
\begin{align}
I_{\ell_1 \ell_2 \ell_3 \ell_4 \ell_5}^{[n]}(r_1, r_2, r_3, r_4), r_5 \equiv \int k^n dk\; j_{\ell_1}(kr_1) j_{\ell_2}(kr_2) j_{\ell_3}(kr_3) j_{\ell_4}(kr_4) j_{\ell_5}(kr_5).
\label{eqn:five}
\end{align}
We use the trick of \S\ref{sec:triple} twice to find
\begin{align}
I_{\ell_1 \ell_2 \ell_3 \ell_4 \ell_5}^{[n]} &= \left(\frac{2}{\pi}\right)^2 \int u'^2 du'\; \int u^2 du\; \int k^n dk\; j_{\ell_1}(kr_1) j_L(ku) j_{L'}(k u')\nonumber\\
& \times \int dk' \; j_{\ell_2}(k' r_2) j_{\ell_3}(k' r_3) j_L(k' u) \; \int dk'' \; j_{\ell_4}(k'' r_4) j_{\ell_5}(k'' r_5) j_{L'}(k'' u').
\label{eqn:five_split}
\end{align}
We now see three triple-SBF integrals that we can perform using the results of \S\ref{sec:triple}. Furthermore, they will involve integrals of at most hypergeometrics, and these integrals, by change of variable, can be seen to result in powers of \textit{e.g.} $u$ and $u'$; hence the final outer integrals can be done in closed form too.

We note that if one wished to perform an integral of six SBFs, one could use our trick of \S\ref{sec:triple} just once; this would give the analog of equation (\ref{eqn:four}) but with four SBFs each in the integrals over $k'$ and $k$ (imagine adding $j_{\ell_5}(k' r_5)$ to the product within the first and $j_{\ell_6}(k r_6)$ to that within the second, for instance). But we have already shown that any integral over four SBFs can be reduced as in equation (\ref{eqn:four}); hence the result of the reduction on the six-SBF integral can be evaluated. Similar logic shows that an arbitrary number of SBFs in the integrand may be dealt with in this fashion.

\section{Concluding Discussion}
\label{sec:concs}
In this work, we have presented a solution to the overlap integral of any pair of spherical Bessel functions with any arbitrary (integer and positive semi-definite) power-law weight. This work has applications in cosmology as more fully outlined in \S\ref{sec:intro}. We have also shown how the ability to do any such double integral enables performing \textbf{any} triple-sBF integral, exploiting a technique originally due to \cite{mehrem_4} but which requires the ability to do an arbitrary double-sBF integral to enable an arbitrary (positive semi-definite integer power-law weight) triple one. Once any triple sBF integral can be done, further use of this Mehrem technique enables doing any quadruple, quintuple, \textit{etc.}, as we outline. An intriguing direction of future work would be to extend these approaches to integrals with negative power-law weights, \textit{i.e.} $k^{-1}$, \textit{etc.} We believe this could be done beginning from the integrals in the present work using a technique of \cite{tomas14}, wherein one performs ``parametric integration'', to lower the power-law. Also of interest would be further investigating the reason for the singularity of some of these overlap integrals but not others; complementary perspectives may come from considering contour integration as an approach, as in \cite{chen04}, and also from analyzing the asymptotic behavior of the phase (\textit{e.g.} as in equation (\ref{eqn:large_asymp})) at large argument [Cahn \& Slepian in prep.]. 

 \section*{Acknowledgments}
We thank J. Chellino,  F. Kamalinejad, A. Krolewski, O.H.E. Philcox, M. Simonovic, and M. White for many useful conversations, and especially R.N. Cahn for much insight into singular sBF integrals over several years, and J. Hou for help in the final stages of various reductions.

\bibliographystyle{unsrt}
\bibliography{refs}

\end{document}